\documentclass[fleqn,10pt]{IEEEtran}
\pdfoutput=1

\usepackage{graphicx}
\usepackage{amssymb}
\usepackage{amsmath}
\usepackage{amsthm}

\usepackage{epstopdf}
\usepackage{enumerate}
\usepackage{algorithmic}
\usepackage{algorithm}
\usepackage{cite}

\DeclareMathOperator{\real}{real}

\theoremstyle{definition}
\newtheorem{defn}{Definition}
\newtheorem{lem}{Lemma}
\newtheorem{prop}{Proposition}
\newtheorem{corollary}{Corollary}

\newcommand{\spa}{8pt}
\title{Sparsity Within and Across Overlapping Groups}
\author{\.Ilker Bayram  \vspace*{-0.5cm}
\\  \thanks{\.{I}. Bayram is with the Dept. of Electronics and Communications Eng., Istanbul Technical University, Istanbul, Turkey.  E-mail : ibayram@itu.edu.tr.
}
}
\date{}

\begin{document}

\maketitle

\begin{abstract}
Recently, penalties promoting signals that are sparse within and across groups  have been proposed. In this letter, we propose a generalization that allows to encode more intricate dependencies within groups. However, this complicates the realization of the threshold function associated with the penalty, which hinders the use of the penalty in energy minimization. We discuss how to sidestep this problem, and demonstrate the use of the modified penalty in an energy minimization formulation for an inverse problem.
\end{abstract}
\vspace*{-0.2cm}
\section{Introduction} \label{sec:intro} 
Sparsity has played a major role in signal processing in the last two decades. However, for many natural signals, plain sparsity falls short of capturing the intrinsic characteristics of the signal of interest. In recent work \cite{bay17px}, we addressed a specific form of sparsity,  useful for signals that are composed of a few number of groups where within each group, only a few coefficients are active. We called this more intricate form of sparsity as `sparsity within and across groups' (SWAG)  -- this characteristic was also referred to as elitist-Lasso \cite{kow09p303,kow09p251}, or exclusive-lasso \cite{zho10AISTATS}, previously in the literature. In this paper, we propose a modification of this penalty that introduces further flexibility in the definition of the groups. We also describe an algorithm to demonstrate how  the proposed penalty can be utilized in simple inverse problem settings.

\subsection{The SWAG Penalty and Threshold Function}
The SWAG penalty in \cite{bay17px} is a group-based penalty. Suppose we are given a collection of variables $x = \{ x_1, x_2,  \ldots, x_N\}$. On this collection, we first form a partition $x^1$, $x^2$, \ldots, $x^k$, where each $x^m$ is a collection of distinct variables from $x$, referred to as a group. Then, the SWAG penalty is defined as
\begin{equation}\label{eqn:P}
P(x) = \|x\|_1 + \frac{\gamma}{2}\,\sum_m \, \sum_{\substack{i,j \\ i \neq j}}  |x^m_i\,x^m_j|.
\end{equation}
The associated threshold function, or proximity operator \cite{Bauschke,combettes_chp} is defined as,
\begin{equation}
T(z) = \arg \min_{x\in \mathbb{C}^n} \frac{1}{2} \|z - x \|_2^2 + \lambda\,P(x). \label{eqn:T}
\end{equation}
$T(z)$ is well-defined if $\lambda \gamma < 1$. $T(z)$ can be computed with a finite terminating procedure \cite{bay17px} and it is group-separable. 

A shortcoming of this penalty is the requirement that the groups be non-overlapping. This constraint is driven primarily by the desire to obtain a realizable threshold function. When the groups share variables, the threshold function is no longer group-separable. In that case, one way to realize the global threshold is to split variables and employ group-separable penalties iteratively in a splitting scheme such as Douglas-Rachford \cite{com08p014,combettes_chp} or ADMM \cite{boy11p1}. Other than the increase in the number of variables, such an approach may not be feasible because some formulations may require to compute infinite iterations within iterations --  a procedure not realizable in principle. 
\begin{figure}
\centering
\includegraphics[scale=1]{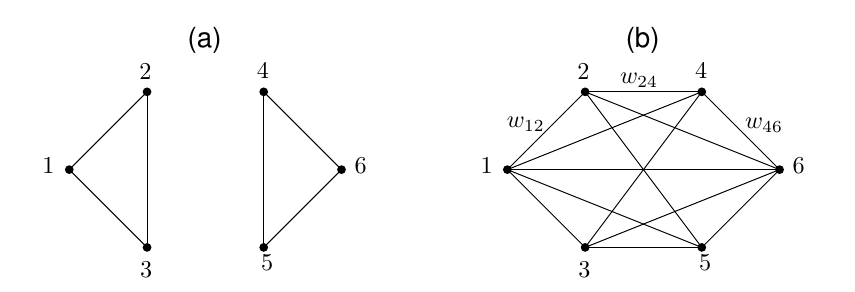}
\caption{A visual description of the group structure using graphs. Each variable is represented by a node. (a) This graph represents the partition $x^{1} = \{x_1, x_2, x_3\}$, $x^{2} = \{x_4,x_5,x_6\}$. 
(b) The proposed generalization employs a weighted complete graph. \label{fig:graphs1}}
\end{figure}

\subsection{The Proposed Penalty}
In order to describe the proposed penalty, we will use graphs as a visual aid. Consider the collection of variables, $x = \{x_1, x_2, \ldots, x_6\}$. Suppose we partition $x$ as $x^{1} = \{x_1, x_2, x_3\}$, $x^2 = \{x_4,x_5,x_6\}$. This partition is represented by the graph in Fig.~\ref{fig:graphs1}a. Notice that each group leads to a complete graph. Since the groups do not share variables, there are two disjoint complete graphs. 
The generalization we propose in this letter is to use a complete weighted graph, as shown in Fig.~\ref{fig:graphs1}c. The modified penalty on $\mathbb{C}^n$ is then defined as
\begin{equation}
P_{W}(x) = \|x\|_1 + \frac{1}{2}\,\sum_{i, j} w_{ij}\,|x_i\,x_j|, \label{eqn:Pw}
\end{equation}
Notice that for a specific choice of $W$, we can recover the penalty \eqref{eqn:P}. Therefore, $P_W$ is a generalization of $P$ in \eqref{eqn:P}.

The associated threshold function is defined similarly as,
\begin{multline}
T_{\lambda,W}(z) = \arg \min_{x\in \mathbb{C}^n} \Bigl\{ D_{\lambda,W}(x;z) \\
= \frac{1}{2} \|z - x \|_2^2 + \lambda\,P_{W}(x) \Bigr\}  \label{eqn:TW}.
\end{multline}
We show in the following section that $T_{\lambda,W}(z)$ is well-defined, provided that $\lambda$ satisfies an upper bound determined by the weight matrix $W$.

\subsection{Related Work and Contribution}
The sparsity characteristic sought in this letter is different than that sought in many papers using group-based penalty functions. Specifically, \cite{yua06p49,kow09p303,jac09ICML,bay11ICASSP,che14p476,sim13p231} aim to promote signals that can be represented with a few groups, where within groups, the coefficients are less stringently penalized. In contrast, \cite{kow09p303,kow09p251,zho10AISTATS,bay17px} aim a similar characteristic as the proposed penalty. In this collection, the SWAG penalty \cite{bay17px}, which the proposed penalty aims to modify, separates from the rest in that it is a non-convex penalty. In \cite{bay17px}, it was argued that this property reduces the bias in the non-zero estimates produced by the threshold function (see also \cite{che14p464,sel14p078} for related discussions). For a more detailed comparison between the SWAG penalty and the penalties in \cite{kow09p303,kow09p251,zho10AISTATS}, we refer to \cite{bay17px}.

The proposed modification to the SWAG penalty aims to introduce further flexibility in forming the groups. First, groups are allowed to overlap. Second, while the original SWAG penalty in \cite{bay17px} uses constant weights within each group, the modified penalty allows the weights within a group to vary. These in turn allow to achieve a more localized and translation-invariant behavior, which is of interest for processing time-domain signals. However, these modifications come at an expense. While it is possible to realize the SWAG threshold function with a finite terminating procedure \cite{bay17px}, such a procedure is not available for the proposed penalty. Therefore, forward-backward splitting type algorithms that might utilize $T_{\lambda,W}$ \cite{combettes_chp,com05p168,fig03p906,bay16p597} are not readily applicable for the proposed penalty. We propose instead a descent algorithm for a generic formulation that employs the proposed penalty. This algorithm is specific to the proposed penalty, and makes use of the quadratic nature of the penalty. Therefore, it has not appeared elsewhere in the literature as far as we are aware.

\subsection*{Notation}
Throughout the manuscript, we take $W$ to be the real symmetric non-negative matrix with entries $W_{i,j} = w_{i,j}$. We assume that the diagonal of $W$ is zero.
For $x\in \mathbb{C}^n$, $x_i$ denotes the $i^{\text{th}}$ component of  $x$, and $|x|$ denotes the vector consisting of the magnitudes of $x_i$'s. We therefore write
\begin{equation}
\sum_{i, j} w_{i,j}\,|x_i\,x_j| = |x|^T\,W\,|x|,
\end{equation}
$\mathbf{1}$ denotes a vector of ones. For non-zero $z\in\mathbb{C}^n $, $e^{j \angle z}$ denotes a unit vector vector in the direction of $z$. For two vectors $x$, $z$ in $\mathbb{C}^n$,  the vector obtained by element-wise multiplication is denoted as $x\,z$.  $\mathbb{R}^n_+$ denotes the non-negative orthant of $\mathbb{R}^n$. 

Finally, $\mathbb{C}^n$ appears as a domain for some functions (including the proposed penalty) in the letter. For inner products and gradients, we interpret $\mathbb{C}^n$ as $\mathbb{R}^{2n}$. Thus, on $\mathbb{C}^n$, we use the inner product $\langle x, y \rangle =  \sum_i \real( x_i y_i^* )$.

\subsection*{Outline}
In Section~\ref{sec:weakconvex}, we derive a condition on $\lambda$ and $W$ which ensures that $P_W$ is weakly-convex and the threshold function $T_W$ is well-defined. Following this, we discuss in Section~\ref{sec:descent} how to construct a descent algorithm when this penalty is used in a simple energy minimization formulation. We demonstrate the utility of the proposed penalty and the minimization algorithm in Section~\ref{sec:Demo}. Section~\ref{sec:conc} contains concluding remarks.

\section{Weak Convexity of the Proposed Penalty}\label{sec:weakconvex}
In this section, we study the proposed penalty function and show that it is weakly convex \cite{via83p231}. 
\begin{defn}
A function $g$ is said to be $\alpha$-weakly convex if 
\begin{equation}
\frac{\alpha}{2}\,\|x\|_2^2 + g(x)
\end{equation}
is convex.
\end{defn}
Our interest in showing the weak convexity of the proposed penalty stems from available schemes such as \cite{sel14p078,bay16p597} that make use of the weak convexity of the penalties. However, as a byproduct of this discussion, we also obtain that $T_{\lambda,W}$ is well-defined provided $\lambda$ is small enough. To see this, observe that $P_{W}$ is $1/\lambda$-weakly convex if and only if $D_{\lambda, W}(\cdot;z)$ (see \eqref{eqn:TW}) is convex. In particular, if $D_{\lambda, W}(\cdot;z)$ is strictly convex, it has a unique minimizer which is in fact $T_{\lambda,W}(z)$. We remark however that $T_{\lambda,W}$ is well-defined for an extended range of $\lambda$ values than those implied by the main result of this section, Prop.~\ref{prop:convex}. We come back to this issue in in Prop.~\ref{prop:well}.

Let us now discuss when $D_{\lambda, W}(\cdot;z)$ is strictly convex.
Notice that, the strict convexity of $|x|^T\,(I + \lambda\, W)\,|x|$ implies the strict convexity of $D_{\lambda, W}(x;z)$ with respect to $x$. We remark that due to the absolute values surrounding $x$, positive definity of $I + \lambda\, W$ does not automatically imply the desired convexity\footnote{Interestingly, in addition to being positive definite,  $I + \lambda\, W$ is also non-negative, and this allows the application of the Perron-Frobenius theorem \cite{Horn} in this context. Therefore, the largest eigenvalue of $I + \lambda\, W$ is unique and the corresponding eigenvector is non-negative. But even this does not appear to imply the desired convexity.}. However, if $I + \lambda\, W$ admits a decomposition of the form
\begin{equation}
I + \lambda\, W = R^T\,R, \text{ with } R_{i,j}\geq 0, \text{ for all } i, j,
\end{equation}
then it can be shown that $D_{\lambda,W}$ is convex. To see this, observe that
\begin{equation}\label{eqn:expand}
|x|^T\,R^T\,R\,|x| = \sum_i \biggl( \sum_j r_{ij}|x_j| \biggr)^2.
\end{equation}
Since $r_{ij}\geq 0$ for all $i,j$, the term enclosed in parentheses in \eqref{eqn:expand} is convex for all $i$, because it is the composition of an increasing function on the positive axis, namely $(\cdot)^2$, and a non-negative convex function, namely a weighted $\ell_1$ norm.

Matrices that admit a decomposition as in \eqref{eqn:expand} are called completely positive \cite{ber88p57}. Unfortunately, checking whether an arbitrary psd matrix is completely positive or not is not a trivial task when the size of the matrix exceeds $4\times4$ \cite{gray80p119,ber88p57}. However, it is relatively simple to find an upper bound for $\lambda$ so that $ I + \lambda W$ is completely positive \cite{ber88p57}. 
\begin{prop}\label{prop:convex}
If
\begin{equation}\label{eqn:constraint}
\lambda\, \Bigl( \max_{i}\,\sum_{j\neq i} w_{i,j} \Bigr) < 1,
\end{equation}
 then $D_{\lambda, W}(\cdot;z)$  is strictly convex.
\begin{proof}
See Appendix~\ref{proof:propconvex}.
\end{proof}
\end{prop}
This proposition implies that $T_{\lambda,W}$ is well-defined if $\lambda$ satisfies \eqref{eqn:constraint}. However, even though $T_{\lambda,W}$ is well-defined, it is not easy to evaluate numerically.  In the sequel, we discuss how to construct descent algorithms for $P_W$. 

\section{Descent Algorithms}\label{sec:descent}

In this section, we derive a descent algorithm for a problem of the form
\begin{equation}\label{eqn:cost}
\min_x\, \bigl\{ C(x) = f(x) + \lambda\,P_{W}(x) \bigr\},
\end{equation}
where $f(\cdot) : \mathbb{C}^n \to \mathbb{R}$ is a convex function. Viewing $\mathbb{C}^n$ as $\mathbb{R}^{2n}$, we also assume that $f$ is Fr\'{e}chet-differentiable \cite{Bauschke,Ortega}, and its Fr\'{e}chet-derivative, $\nabla f$, is Lipschitz-continuous with parameter $L$, i.e.,
\begin{equation}
\| \nabla f(x) - \nabla f(y) \|_2 \leq L \|x - y \|_2, \text{ for all } x, y.
\end{equation}

We will derive the algorithm based on the majorization-minimization scheme \cite{hun04p30, fig07p980}. Specifically, we will discuss how to update the $k^{\text{th}}$ iterate $x^k$ so that $C(x^{k+1}) \leq C(x^k)$.
\begin{defn}\label{defn:majorizer}
A function $g : \mathbb{C}^n \to \mathbb{R}$ is said to be a \emph{majorizer for $h:\mathbb{C}^n \to \mathbb{R}$ at $x^*$} if
\begin{enumerate}[(i)]
\item $h(x^*) = g(x^*)$,
\item $h(x) \leq g(x)$ for all $x \in \mathbb{C}^n$.
\end{enumerate}  
\end{defn}
We first provide a majorizer for $f$. Although that is more or less well-known, we include a short discussion for the sake of completeness.

\subsection{Majorizing `$f(\cdot)$'}
Thanks to the properties of $f$, we can readily obtain a majorizer as follows.
\begin{prop}\label{lem:majorizef}
Suppose $f$ is convex and its Fr\'{e}chet derivative $\nabla f$ is Lipschitz continuous with parameter $L$. If $\alpha \leq 1/L$,  then
\begin{multline}\label{eqn:Ck}
C^k(x) =  \frac{1}{2\alpha} \left\|x -  \left(x^k - \alpha \nabla f(x^k) \right) \right\|_2^2 + \lambda\,P_w(x) \\
+ \Bigl[ f(x^k) - \frac{\alpha}{2}\,\| \nabla f(x^k) \|_2^2 \Bigr].
\end{multline}
is a majorizer for $C(\cdot)$ at $x^k$.
\begin{proof}
See Appendix~\ref{proof:majorizef}.
\end{proof}
\end{prop}
Notice that the term inside the square brackets in \eqref{eqn:Ck} is constant with respect to $x$ and does not play a role in the subsequent minimization.
By the two properties of a majorizer, it follows that if we set $x^{k+1}$ to be a minimizer of $C^k(x)$, then
$C(x^{k+1}) \leq C(x^k)$.
However, to obtain $x^{k+1}$, we essentially need to solve the problem in \eqref{eqn:TW}, for which a numerical procedure is not readily available. Nevertheless, thanks to the two properties of $C^k$ listed above, if we find $x^*$ that achieves $C^k(x^*)\leq C^k(x^k)$, we will have
$C(x^*) \leq C(x^k)$.
Thus, it suffices to perform descent on $C^k$. In the following, we show that this can be achieved with a simple update rule.

\subsection{Majorizing the Proposed Penalty}
In view of the foregoing discussion, our goal is to find  some $\hat{x}$ such that $C^k(\hat{x}) \leq C^k(x^k)$. This condition is equivalent to 
\begin{equation}
D_{\beta,W}(\hat{x};z) \leq D_{\beta,W}(x^k ; z)
\end{equation}
for $z = x^k - \alpha \nabla f(x^k)$, and $\beta = \alpha\,\lambda$.

Observe now that
\begin{equation}\label{eqn:observe}
D_{\beta,W}\bigl( |x|;|z| \bigr) = D_{\beta,W}\bigl(|x|\,e^{j\angle z} ; z \bigr) \leq D_{\beta,W}\bigl( x ; z \bigr).
\end{equation}
for all $x$, $z$ in $\mathbb{C}^n$.
This suggests that, instead of minimizing $D_{\beta,W}(x;z)$, we can consider a minimization problem as
\begin{equation}\label{eqn:probconstrained}
\min_{x \in \mathbb{R}^n_+}\, D_{\beta,W}(x;|z|).
\end{equation}
On $\mathbb{R}^n_+$, $D_{\beta,W}(\cdot ; |z|)$ is simply a quadratic function. 
This has the following consequence.
\begin{prop}\label{prop:well}
If $I +  \beta\,W$ is positive definite, then $D_{\beta,W}(\cdot ; z)$ has a unique minimizer.
\begin{proof}
See Appendix~\ref{proof:propwell}.
\end{proof}
\end{prop}
Notice that this extends the range implied by Prop.~\ref{prop:convex} over which $T_{\lambda,W}$ is well-defined. However, it does not imply the strict convexity of $D_{\beta,W}(\cdot;z)$, as Prop.~\ref{prop:convex} does.

The problem in \eqref{eqn:probconstrained} is a constrained convex minimization problem. Thus, descent can be achieved by applying any finite number of iterations of the projected gradient algorithm \cite{gol64p709}. This observation leads to the following result.
\begin{prop}\label{prop:descent}
Suppose $ I + \beta\,W$ is positive semi-definite with spectral norm $\sigma$; $f$ is convex, and its Fr\'{e}het derivative $\nabla f$ is Lipszhitz continuous with parameter $L$; and $\alpha \leq 1/L$. Let $z = x^k - \alpha\,\nabla f(x^k)$, and $P_+(\cdot)$ denote the projection operator onto $\mathbb{R}^n_+$. Also let $S:\mathbb{R}^n \to \mathbb{R}^n_+$ denote the operator that maps $x\in \mathbb{R}_+$, to $\hat{x} \in \mathbb{R}^n_+$ where,
\begin{equation}
\hat{x} = P_+ \left( x - \eta \left[ \Bigl( I + \beta\,W  \Bigr) x  + \beta \mathbf{1} - |z| \right]  \right).
\end{equation}
Finally, let $S^m$ denote $S$ iterated $m$ times.
If $\eta \leq 2 / \sigma$, then for any $m\geq 1$, and $C$ as in \eqref{eqn:cost}, we have
\begin{equation}\label{eqn:descentguarantee}
C\bigl( S^m(|x^k|)\,e^{j \angle z} \bigr) \leq C\bigl( x^k \bigr).
\end{equation}
Further, if equality holds in  \eqref{eqn:descentguarantee}, then,
\begin{enumerate}[(i)]
\item $x^k$ is a stationary point of $C(\cdot)$, i.e., 0 is in the proximal subdifferential \cite{Clarke2} of $C(\cdot)$,
\item $x^k = S^m(|x^k|)\,e^{j \angle z}$, i.e., $x^k$ is a fixed point of the iterations.
\end{enumerate}
\begin{proof}
See Appendix~\ref{proof:descent}
\end{proof}
\end{prop}
In view of  Prop.~\ref{prop:descent}, Algorithm~\ref{algo:Descent} achieves descent for \eqref{eqn:cost}.

\begin{algorithm}[h!]{\caption{A Descent Algorithm for \eqref{eqn:cost}}\label{algo:Descent}}
\begin{algorithmic}[1]
\REQUIRE $L$, Lipschitz const. of $\nabla f$; $K$, number of inner iterations ;
$\lambda$, weight of $P_W$
\STATE Set $\alpha < 1/L$, $\beta \gets \alpha\,\lambda$,  $\eta < 2 / \sigma\bigl(I +  \beta\,W\bigr)$, initialize $x$
\REPEAT
\STATE $z \gets x - \alpha \nabla f(x)$\label{line:nabla}
\STATE $x \gets |x|$ 
\FOR{$K$ iterations}
\STATE $x \gets P_+ \Bigl( x - \eta \bigl[(\,I + \beta\,W)x + \beta \mathbf{1} - |z| \bigr] \Bigr)$
\ENDFOR
\STATE $x\gets x \, e^{j\,\angle z }$
\UNTIL{convergence}
\end{algorithmic}
\end{algorithm}

\section{Demonstration of the Proposed Penalty/Algorithm}\label{sec:Demo}
In this section, we demonstrate the utility of the proposed penalty and the algorithm on a dereverberation experiment. Our purpose is to show the differences of the proposed penalty compared to the SWAG penalty. 

\begin{figure}
\centering
\includegraphics[scale=1]{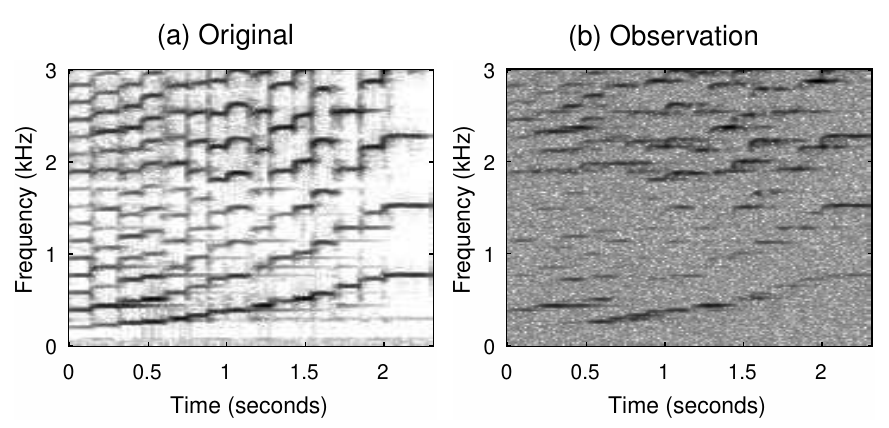}
\caption{Spectrograms of (a) the original signal, and (b) the reverberant and noisy observation used in the experiment.\label{fig:spectrogram}}
\end{figure}

The clean signal is comprised of a violin playing a chromatic scale, sampled at 16 KHz. The spectrogram of the clean signal is shown in Fig.~\ref{fig:spectrogram}a. The computations are carried out in the STFT domain, so that the effects of the penalty are easier to observe. The reverberant spectrogram is obtained by convolving each STFT band with a filter obtained from an impulse response \cite{rei02p730}. This amounts to applying a linear operator, say $H$, to the clean STFT coefficients. We then add circular complex valued Gaussian noise to the reverberant signal's STFT coefficients so that the observation SNR is 5~dB (see Fig.~\ref{fig:spectrogram}b).

We consider a reconstruction formulation as
\begin{equation}
\min_x  \underbrace{\frac{1}{2} \| y - H\,x \|_2^2}_{f(x)} + P(x),
\end{equation}
where $P(x)$ is the penalty term, which is one of $\ell_1$ norm, the SWAG penalty or the proposed penalty.

\begin{figure}
\centering
\includegraphics[scale=1]{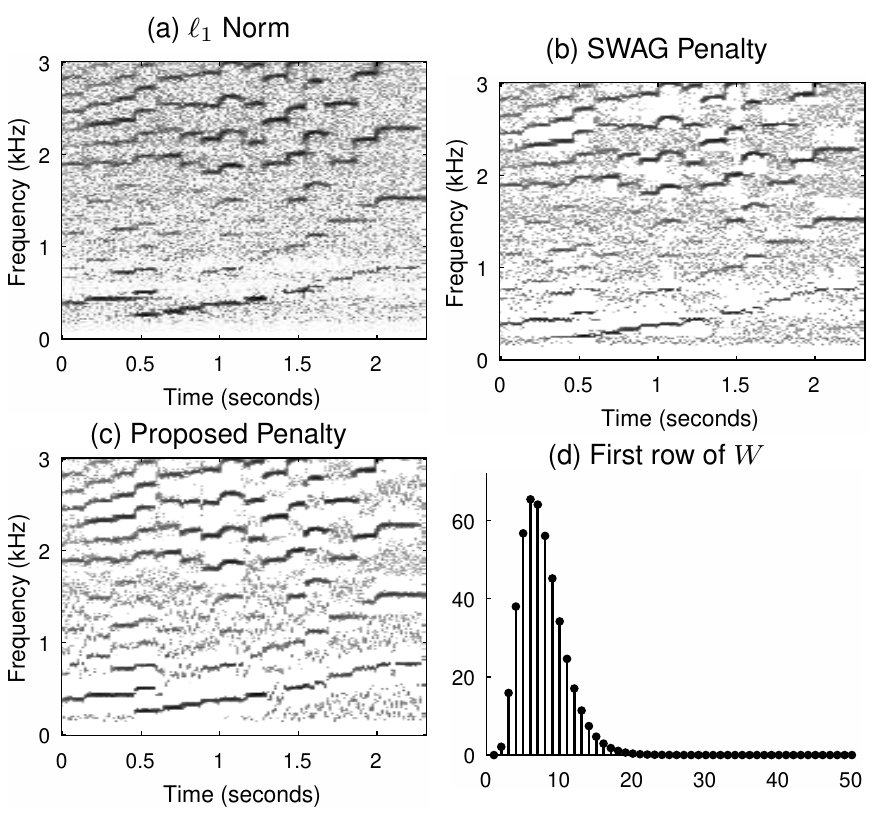}
\caption{Spectrograms of the dereverbed signals using (a) the $\ell_1$ norm, (b) SWAG penalty from \cite{bay17px}, (c) proposed penalty. (d) The first 50 coefficients from the first row of the Toeplitz weight matrix $W$, used for defining $P_W(\cdot)$.\label{fig:reconstruction}}
\end{figure}

We select the weight of the $\ell_1$ norm, with a sweep search so as to maximize the SNR.
For both SWAG and the proposed penalty, we set $\lambda$ to be a quarter of the weight used for the $\ell_1$ penalty. For these penalties, to form a group, we take a perpendicular (referring to Fig.~\ref{fig:spectrogram},~\ref{fig:reconstruction}) slice (i.e., along the frequency axis) in the STFT domain. For SWAG, we partition this slice into groups of size 15 and set $\gamma$ to 40. For the proposed penalty, we set the weight matrix as a symmetric Toeplitz matrix of size $960\times 960$\footnote{960 is the total number of frequency bands.}, where the first column is as shown in Fig.~\ref{fig:reconstruction}d. The sequence is non-negative and sums to 900, so that $\sigma(W) \leq  900$. A Toeplitz $W$ with non-zeros close to the main diagonal leads to a localized effect in terms of penalizing coefficients and achieves translation-invariance.

We remark that, for $f$, the Lipschitz parameter $L$ is equal to the largest eigenvalue of $H^*\,H$, which is approximated numerically. The parameter $\alpha$ is set to $0.9 / L$ to guarantee convergence. We also use the upperbound stated in the text for the spectral norm of $W$, $\sigma_W$, and set $\eta = 1.9 / (1 + \beta\,\sigma_W)$.

The reconstructions obtained by using the three different regularizers are shown in Fig.~\ref{fig:reconstruction}. The output SNRs are, 7.75 for the $\ell_1$ norm, 5.76~dB for the SWAG penalty and 7.80~dB for the proposed penalty. We remark that our purpose here is not to compare the methods based on the output SNR, but to demonstrate the different characteristics of the reconstructions via the spectrograms.

For the $\ell_1$ norm, if we use a higher threshold, that leads to the suppression of the weaker harmonics, along with noise. SWAG and the proposed penalty avoid this dilemma by suppressing noise around the strong harmonics and retaining the weaker harmonics, even if they are surrounded by noisy coefficients. Overall, this still leads to an improvement in terms of SNR, at least for the proposed penalty. This aside, the spectrograms obtained with SWAG and the proposed penalty show some differences. SWAG uses non-overlapping groups. Also, since the boundaries of the groups are not selected with respect to the positions of the harmonics, we observe that the cleared portions surrounding the harmonics are not centered around the harmonics. In contrast, thanks to the Toeplitz nature of $W$, the proposed penalty essentially employs maximally overlapping groups. This leads to a reconstruction where the harmonics lie at the center of an otherwise suppressed area. This is especially easier to observe in the harmonics occurring after 2 sec.

\section{Outlook}\label{sec:conc}
The proposed penalty allows enhanced flexibility in forming the groups and has the potential to enhance the reconstruction performance that can be obtained with the SWAG penalty. In \cite{bay17px}, we have also shown that the SWAG penalty can be used in combination within a higher-level group-forming strategy to obtain `hybrid penalties'. Such modifications are also feasible for the penalty proposed in this paper. 

Another aspect of interest is the selection of the weight matrix $W$. While the proposed penalty offers flexibility in the choice of the groups via the introduction of $W$, it is not obvious what the optimal weights are. One possible approach for the selection of $W$ might be to \emph{learn} it from data. We hope to investigate this problem in  future work.

\appendices
\section{Proof of Prop~\ref{prop:convex}}\label{proof:propconvex}
First, observe that
\begin{equation}
\frac{\beta}{2} (|x_i|^2 + |x_j|^2) + \beta \, | x_i\,x_j| = \frac{\beta}{2}(|x_i| + |x_j|)^2
\end{equation}
is convex for any $\beta \geq 0$. Therefore,
\begin{equation}
\lambda\,\sum_i \sum_{\substack{j \\ j \neq i}} \left( w_{ij}\,|x_i\,x_j| + \frac{w_{ij}}{2} (|x_i| + |x_j|)^2 \right)
\end{equation}
is convex. But this function can be expressed as $|x|^T\,(D + \lambda\,W )|x|$, where $D$ is a diagonal matrix with 
\begin{equation}
D_i  = \lambda \,\sum_{\substack{j \\ j\neq i}} w_{ij} < 1.
\end{equation}
Thus,
\begin{equation}
|x|^T\,( I + \lambda\, W)\,|x| = x^H (I - D )x + |x|^T\,(D + \lambda\,W )|x|,
\end{equation}
is a sum of a strictly convex and a convex function, so it is strictly convex.

Observe now that $D_{\lambda, W}$ can be expressed as
\begin{equation}
\Bigl[ \frac{1}{2}\|z\|_2^2 - 2 \langle z,x \rangle + \lambda\,\|x\|_1 \Bigr] + \Bigl\{ |x|^T\,( I + \lambda\, W) |x| \Bigr\}.
\end{equation}
The term inside the square brackets is convex with respect to $x$, and the term inside the curly brackets was shown to be strictly convex. Thus follows the claim.

\section{Proof of Prop.~\ref{lem:majorizef}}\label{proof:majorizef}
Let us start by recalling a lemma from convex analysis. 
\begin{lem}\label{lem:bound}
Suppose $g:\mathbb{R}^n \to \mathbb{R}$ is a convex, Fr\'{e}chet differentiable function whose Fr\'{e}chet derivative $\nabla f$ is Lipschitz continuous with parameter $L$. Then,
\begin{equation}\label{eqn:bound}
g(x) \leq g(y) + \langle \nabla g(y), x-y \rangle + \frac{L}{2}\, \| x - y \|_2^2, \text{ for all } x, y.
\end{equation}
\begin{proof}
See for instance Cor.18.14, (i)$\Rightarrow$(iv) in \cite{Bauschke}.
\end{proof}
\end{lem}
Using this lemma, we obtain the following corollary after some algebra.
\begin{corollary}\label{cor:bound}
Suppose $g$ is as in Lemma~\ref{lem:bound}. Then,
\begin{equation}
\frac{L}{2}\, \Bigl\| x - \Bigl(y - \frac{1}{L} \nabla g(y)\Bigr) \Bigr\|_2^2 + g(y) - \frac{1}{2L} \| \nabla g(y) \|_2^2.
\end{equation}
is a majorizer for $g$ at $y$.
\end{corollary}
Observe that if $\nabla f$ is $L$-Lipschitz continuous then it is also Lipschitz continuous with $L' \geq L$. Therefore, applying Lemma~\ref{cor:bound} to $f$ in $C$  with $y = x^k$, we find that $C^k$ is a majorizer for $C$ at $x^k$.

\section{Proof of Prop.~\ref{prop:well}}\label{proof:propwell}
Notice that the problem in \eqref{eqn:probconstrained} can be written as
\begin{equation}\label{eqn:probmodified} 
\min_{x \in \mathbb{R}^n_+}\,\frac{1}{2} \, x^T\,(I + \beta\, W)\,x - \bigl\langle |z| +  \beta \mathbf{1}, \, x \bigr\rangle.
\end{equation}
Suppose $I + \beta W$ is positive definite. Then, the function to be minimized in \eqref{eqn:probmodified} is strictly convex. Since the constraint set is convex, it follows that \eqref{eqn:probmodified} has a unique solution. Let $x^*\in \mathbb{R}^n_+$ denote the unique solution of \eqref{eqn:probmodified}. We claim that $D_{\beta,W}(x^*\,e^{j\angle z}, z) < D_{\beta,W}(x,z)$ for any $x \neq x^*\,e^{j\angle z}$. Suppose 
\begin{equation}\label{eqn:ineqpf}
D_{\beta,W}(x^*\,e^{j\angle z}, z) \geq D_{\beta,W}(x,z).
\end{equation} 
Invoking \eqref{eqn:observe}, we have $D_{\beta,W}(x^*\,e^{j\angle z}, z) \geq D_{\beta,W}(|x|\,e^{j\angle z},z)$. But again by \eqref{eqn:observe}, this implies $D_{\beta,W}(x^*, |z|) \geq D_{\beta,W}(|x|,|z|)$. By the uniqueness of the solution of \eqref{eqn:probmodified}, we conclude that $x^* =  |x|$.
Observe now that if $e^{j \angle x} \neq e^{j\angle z}$, then $D_{\beta,W}(|x| e^{j\angle z}, z) < D_{\beta,W}(x,z)$. Consequently,
\begin{equation}
D_{\beta,W}(x^* e^{j\angle z}, z) = D_{\beta,W}(|x| e^{j\angle z}, z) < D_{\beta,W}(x,z),
\end{equation}
contradicting \eqref{eqn:ineqpf}. Thus, $e^{j \angle x} = e^{j\angle z}$, and $x = |x| \, e^{j\angle z} = x^*\,e^{j \angle z}$, proving the claim that the minimizer of $D_{\beta,W}(\cdot,z)$ is unique.

\section{Proof of Prop.~\ref{prop:descent}}\label{proof:descent}
Let us first consider the descent property. We already noted in the text that $C^k(x^*) \leq C^k(x^k)$ is equivalent to $D_{\beta,W}(x^*;z) \leq D_{\beta,W}(x^k ; z)$. So it suffices to show the validity of this inequality for $x^* =  S^m(|x^k|)\,e^{j \angle z}$.

In fact, thanks to \eqref{eqn:observe}, it is sufficient to show 
\begin{equation}\label{eqn:suf}
 D_{\beta,W}(S^m(|x^k|);|z|) \leq D_{\beta,W}(|x^k|;|z|)
 \end{equation}
 because if it is valid, then
\begin{subequations}
\begin{align}
D_{\beta,W}(x^*;z\bigr) &= D_{\beta,W}\bigl(S^m(|x^k|);|z|)\\
& \leq D_{\beta,W}(|x^k|;|z|) \\
& \leq D_{\beta,W}(x^k;z).
\end{align}
\end{subequations}

An application of $S$ amounts to one iteration of the projected gradient algorithm on \eqref{eqn:probconstrained}. For the sake of completeness, we also show this implies the claimed descent property. 
Now, let $x\in \mathbb{R}^n_+$.
First, observe that
\begin{equation}
 \nabla D_{\beta,W}(x, |z|) 
= ( I + \beta\,W )\, x  + \beta \mathbf{1} - |z|,
\end{equation}
where differentiation is performed with respect to the first variable.
By the properties of the projection operator onto $\mathbb{R}^n_+$ \cite{HiriartFund}, along with $x \in \mathbb{R}^n_+$, $S(x) \in \mathbb{R}^n_+$, we have
\begin{equation}
\bigl\langle x - S(x), x - \eta  \nabla D_{\beta,W}(x, |z|)  - S(x)  \bigr\rangle \leq 0.
\end{equation}
Rearranging, and invoking Lemma~\ref{lem:bound}, we have,
\begin{subequations}
\begin{align}
\frac{1}{\eta} &\bigl\| x - S(x) \bigr\|_2^2 \nonumber \\
&\leq \bigl\langle x - S(x), 
 \nabla D_{\beta,W}(x, |z|)  \bigr\rangle \\
 &\leq D_{\beta,W}\bigl(x, |z|\bigr) - D_{\beta,W}\bigl(S(x), |z| \bigr) \\
 & \quad \quad + \frac{\sigma}{2} \bigl\| x - S(x) \bigr\|_2^2
\end{align}
\end{subequations}
Rearranging, we obtain
\begin{multline}\label{eqn:Sm}
D_{\beta,W}\bigl(x, |z|\bigr) - D_{\beta,W}\bigl(S(x), |z| \bigr) \\ \geq \Bigl(\frac{1}{\eta} - \frac{\sigma}{2} \Bigr)\, \bigl\| x - S(x) \bigr\|_2^2.
\end{multline}
From the assumption $\eta \sigma < 2$, it follows that the rhs is non-negative. Repeatedly invoking this inequality, we can thus write
\begin{equation}
D_{\beta,W}\bigl(S^m(x), |z| \bigr)\leq D_{\beta,W}\bigl(x, |z|\bigr).
\end{equation}
Plugging in $x = |x^k|$, \eqref{eqn:suf} follows.

Suppose now equality holds in \eqref{eqn:descentguarantee}. This implies that $D_{\beta,W}(S^m(|x^k|);|z|) = D_{\beta,W}(|x^k|;|z|)$. But by \eqref{eqn:Sm}, this is possible only if $S(|x^k|) = |x^k|$. This in turn implies 
\begin{equation}\label{eqn:opt1}
|x^k| - |z| + \beta\mathbf{1} +  W |x^k| \in - N_+(|x^k|),
\end{equation}
where $N_+(|x^k|)$ is the normal cone  \cite{HiriartFund} of $\mathbb{R}^n_+$ at $|x^k|$. 
\eqref{eqn:opt1}  coincides with the optimality condition for $|x^k|$ for   \eqref{eqn:probconstrained}. It then follows from the train of inequalities in \eqref{eqn:observe} that $x^* = S^m(|x^k|)\,e^{j \angle z}$ minimizes $D_{\beta,W}(\cdot;z)$. 
But $C(x^*) =  C(x^k)$ implies $C^k(x^k) \leq C^k(x^*)$, which is equivalent to $D_{\beta,W}(x^*;z) \leq D_{\beta,W}(x^k;z)$. Therefore, $x^k$ also minimizes $D_{\beta,W}(\cdot;z)$. Thus
\begin{equation}
0 \in x^k - z +  \beta\,\partial P_W(x^k),
\end{equation}
where $\partial P_W(x^k)$ is the proximal subdifferential of $P_W$ \cite{Clarke2}.
Plugging in $z = x^k - \alpha\,\nabla f(x^k)$ and $\beta = \alpha\,\lambda$, we  obtain,
\begin{equation}
0 \in \nabla f(x^k) + \lambda\,\partial P_W(x^k).
\end{equation}
Thus, $x^k$ is a stationary point of $C(\cdot)$, as claimed in (i).

Observe now that if $x^k$ minimizes $D_{\beta,W}(\cdot;z)$, then we must have $e^{j\angle z} = e^{j\angle x^k}$. Therefore,   $x^* = S(|x^k|)|\,e^{j\angle z} = |x^k|\,e^{j\angle x^k} = x^k$, as claimed in (ii).

\end{document}